\documentclass{amsart}
\usepackage{amssymb}
\usepackage[dvips]{epsfig}
\usepackage{graphicx}
\usepackage{color}
\usepackage{amsmath}
\usepackage{amssymb}         
\usepackage{amsfonts}
\usepackage{amsthm}
\usepackage{bbm}

\newcommand{\R}{\mathbb R}
\newcommand{\p}{\partial}

\newcommand{\Z}{\mathbb Z}

\newcommand{\C}{\mathbb C}

\renewcommand{\a}{\alpha}

\newcommand{\tres}{|\!|\!|}

\renewcommand{\u}{\mathfrak{u}}

\newcommand{\ji}{\langle}
\newcommand{\jd}{\rangle}

\newtheorem{theorem}{Theorem}[section]
\newtheorem{lemma}[theorem]{Lemma}

\newtheorem{corollary}[theorem]{Corollary}
\theoremstyle{remark}
\newtheorem{remark}{Remark}[section]
\theoremstyle{definition}

\numberwithin{equation}{section}

\begin{document}
\title[Generalized derivative Schr\"odinger equation]{On a class of solutions to the generalized derivative Schr\"odinger equations II}
\author{F.  Linares}
\address[F. Linares] {IMPA\\ Estrada Dona Castorina 110, Rio de Janeiro 22460-320, RJ Brazil}
\email{linares@impa.br}
\author{G. Ponce}
\address[G. Ponce]{Department  of Mathematics\\
University of California\\
Santa Barbara, CA 93106\\
USA}
\email{ponce@math.ucsb.edu}
\author{G.N. Santos}
\address[G.N. Santos]{Universidade Federal do Piau\'i - UFPI.
Campus Ministro Petr\^onio Portella\\
Teresina 64049-550, PI Brazil}
\email{gleison@ufpi.edu.br}

\thanks{F.L. was partially supported by CNPq and FAPERJ/Brazil}

\begin{abstract}
In this note we shall continue our study on the initial value problem associated for the generalized derivative Schr\"odinger  (gDNLS) equation
\begin{equation*}
\p_tu=i\p_x^2u + \mu\,|u|^{\a}\p_xu, \hskip10pt x,t\in\R,  \hskip5pt 0<\a \le 1\;\; {\rm and}\;\; |\mu|=1.
\end{equation*}

Inspiring by Cazenave-Naumkin's  works we shall establish the local well-posedness for
a class of data of arbitrary size in an appropriate weighted Sobolev space, thus removing the size restriction on the data required in our previous work. The  main new tool in the proof  is the homogeneous
and inhomogeneous versions of the Kato smoothing effect for the linear Schr\"odinger equation with lower order variable coefficients established by Kenig-Ponce-Vega.
\end{abstract}

\maketitle




\section{Introduction}

We shall deal with the initial value problems (IVP) associated to the generalized nonlinear derivative 
Schr\"odinger (gDNLS) equation,
\begin{equation}\label{gdnls}
\begin{cases}
\p_tu=i\,\p_x^2u + \mu\,|u|^{\a}\p_xu, \hskip15pt x,t\in\R, \hskip5pt 0<\a \le 1,\\
u(x,0)=u_0(x),
\end{cases}
\end{equation}
where $u$ is a complex valued function, $\mu\in \C$ with $|\mu|=1$.

The equation in \eqref{gdnls} generalized the well-known derivative nonlinear Schr\"odinger  (DNLS) equation
\begin{equation}\label{dnls-0}
i\p_tu+\p_x^2u + i\p_x(|u|^2 u)=0, \hskip15pt x,t\in\R,
\end{equation}
which arises as a model in plasma physics and optics (\cite{Mio}, \cite{Mijolhus}, \cite{MMW}).  It is also an equation which is  exactly solvable by the inverse scattering technique, see \cite{KN}. 

The IVP associated to \eqref{dnls-0} has been considered
in several publications  (see for instance \cite{Linares}, \cite{Iteam}, \cite{refinement}, \cite{HayashiII}, \cite{HayashiOzawa}, \cite{HayashiOzawa2}, \cite{Ozawa}, \cite{Takaoka},  \cite{Takaoka2}, \cite{Tsutsumi}) where among other qualitative properties local and global well-posedness issues were investigated. In particular,  a global sharp well-posedness result was obtained by  Guo and Wu \cite{GW} in $H^{1/2}(\R)$ for initial data satisfying appropriate restrictions, (see also \cite{W}).

 Recently, via inverse scattering method Jenkins, Liu, Perry and Sulem \cite{JLPS}  established global existence of solutions without any restriction on the size of the data in an appropriate weighted Sobolev space.  For results concerning the initial-periodic-boundary value problem (IPBVP)  we refer to \cite{Grunrock} and \cite{Herr}.

Similar to the DNLS equation the gDNLS admits a two-parameter family of solitary wave solutions given explicitly by 
\begin{equation*}
\psi_{\omega,c}(x,t)=\varphi_{\omega,c}(x-ct) \exp i\Big\{ \omega t+\frac{c}{2}(x-ct)-\frac{1}{\a+2} \int\limits_{-\infty}^{x-ct} \varphi_{\omega,c}^{\a}(y)\,dy\Big\}, \hskip10pt \a>0,
\end{equation*}
where
\begin{equation*}
\varphi_{w,c}(x)=
\begin{cases}
\Big\{\dfrac{(2+\a)(4\omega-c^2)}{4\sqrt{\omega}\cosh(\frac{\a}{2}\sqrt{4\omega-c^2}x\big)-2c}\Big\}^{\frac{1}{\a}},  \hskip20pt \text{if} \hskip10pt \omega>\frac{c^2}{4},\\
\\
(\a+2)^{\frac{1}{\a}} \,c^{\frac{1}{\a}}\,(\frac{\a^2}{4}(cx)^2+1)^{-\frac{1}{\a}}, \hskip7pt \text{if} \hskip10pt \omega=\frac{c^2}{4} \hskip10pt \text{and} \hskip10pt  c>0.
\end{cases}
\end{equation*}

For the study of the stability  and instability  for these solitary wave solutions we refer to \cite{LSS}, \cite{TX}, \cite{Guo}, \cite{Colin}, \cite{KW}, \cite{MTX}, \cite{LW} and references
therein.
 
Regarding the well-posedness of the IVP  \eqref{gdnls} with  $\a>0$, $\a\ne2$, in \cite{Hao} Hao obtained
 local well-posedness  in $H^s(\R)$, for $\a>5$ and $s>\frac12$.  In \cite{Gleison} Santos considered the case of
sufficient small initial data in the case  and $1<\a<2$ and showed
the existence and uniqueness of solution $u=u(x,t)$ with  $$ u\in L^{\infty}((0,T) : H^{3/2}(\R)),\;\,\ji x\jd^{-1}u\in L^{\infty}((0,T):H^{1/2}(\R))$$  and local well-posedness in $H^{1/2}(\R)$ for small data when $\a>2$. In \cite{ MHO}  Hayashi and Ozawa
considered the gDNLS equation  in a bounded interval with a Dirichlet condition and established local results in $H^2$ for $\a\ge1$ and $H^1$ for $\a\ge 2$. In \cite{FHI}   Fukaya, Hayashi and Inui showed global result for initial data in $H^1(\R)$, for any $\a\ge 2$, with initial data satisfying some size restriction.

\vskip.1in

In \cite{Cazenave}  Cazenave and Naumkin studied the IVP associated to semi-linear Schr\"odinger equation,
\begin{equation}\label{nls}
\p_tu=i(\Delta u\pm|u|^{\a} u), \hskip15pt x\in\R^N,\;t\in\R, \hskip5pt \a>0.
\end{equation}
For every $\a>0$ they constructed a class of initial data for which they can prove the existence of a unique local solutions for the IVP \eqref{nls}.  Also, they obtained 
a class of initial data for $\alpha>\frac{2}{N}$ for which there exist global solutions that scatter.

One of the ingredients in their argument is the fact that solutions of the linear problem satisfy
\begin{equation}\label{lower}
\underset{x\in\R^N}{\rm Inf}\; \ji x\jd^m\,| e^{it\Delta} u_0(x)| >0,
\end{equation}
for $t\in[0,T]$ with $T$ sufficiently small whenever the initial data satisfy
\begin{equation}\label{lower-2}
\underset{x\in\R^N}{\rm Inf}\; \ji x\jd^m\,|u_0(x)|\ge \lambda>0.
\end{equation}
This is extended to solutions of the \eqref{nls} with data $\,u_0\,$ satisfying \eqref{lower-2} for $m=m(\a)$ and $u_0\in H^s(\R^N)$ with $s$ sufficiently large.

Stimulated by the ideas in \cite{Cazenave} in our previous paper  \cite{LiPoSa} we established a local well-posedness result for the IVP \eqref{gdnls} for small data in a Sobolev weighted space with data
satisfying \eqref{lower-2}.  
Our main objetive here is to remove the size restriction on the data in our previous paper \cite{LiPoSa}. 
More precisely,  our main result reads as follows.



\begin{theorem}\label{main}
There exists $\,M\in\Z^+$ such that for any $\a\in (0, 1)$, for any $k\ge   m+M+1$, \hskip3pt $k\in\Z^{+}$, with 
$$\,m\equiv \Big[\frac{2}{\alpha}+1\Big] $$
{\rm(}$\,[x]$  the greatest integer less than or equal to $\,x${\rm)}, 
and any data $\,u_0$ satisfying
\begin{equation}
\label{sob}
u_0\in H^s(\R), \hskip15pt s=k+\frac12,
\end{equation}
and
\begin{equation}
\label{newspace}
\ji x\jd^m u_0\in  L^{\infty}(\R),\hskip15pt  \ji x\jd^m \partial_x^{j+1}u_0\in  L^{2}(\R),\;\;\;\;\;\;\;j=0,1,..,M,
\end{equation}

with

\begin{equation}\label{main-1}
\|u_0\|_{s,2}+\| \ji x\jd^m u_0\|_{\infty}+\sum_{j=0}^M\| \ji x\jd^m \partial_x^{j+1}u_0\|_2 \equiv \nu
\end{equation}
and
\begin{equation}\label{main-2}
\underset{x}{\rm Inf}\; \ji x\jd^m\,|u_0(x)|\ge \lambda>0,
\end{equation}
 the IVP  \eqref{gdnls} has a  unique solution $u$ such that
\begin{equation}\label{main-3}
u\in C([0,T] : H^s(\R)),
\end{equation}
\begin{equation}\label{main-3a}
\ji x\jd^m u\in L^{\infty}([0,T] :  L^{\infty}(\R)),
\end{equation}
\begin{equation}\label{main-3b}
\ji x\jd^m\partial_x^{j+1}u\in C([0,T] :L^2(\R)),\;\;\;\;\;\;\;\;j=0,1,...,M,
\end{equation}
and
\begin{equation}\label{main-4}
\sup_{j\in\Z} \|\p_x^{k+1}u\|_{L^2([j,j+1]\times[0,T])}  <\infty,
\end{equation}
with $T=T(\alpha; k; \nu; \lambda)>0$. Moreover, the map data-solution
\begin{equation*}
u_0\mapsto u(\cdot,t)
\end{equation*}
from a neighborhood of the datum $u_0$ in  $H^s(\R)$  intersected with the set in \eqref{newspace}  satisfying \eqref{main-1}-\eqref{main-2} into the class defined by 
\eqref{main-3}-\eqref{main-4}
is locally continuous.
\end{theorem}

\vskip.25in
\begin{remark}

(a) The choice $\,m=[2/\a+1]\,$ can be replaced by $m>1/\a\geq 1$.  However, by fixing  $m=[2/\a+1]$ one greatly abbreviates   some technical details in the proof. 
\vskip.1in
(b) Clearly $\lambda< \nu$. Typically the function 
$\varphi(x)=\dfrac{c_0}{\ji x\jd^m}$ with $c_0\neq 0$ satisfies the hypotheses \eqref{sob}, \eqref{newspace} and \eqref{main-2}. Notice that the traveling wave solutions $\,\psi_{\omega,c}(x,t)\,$ described above are not in this class. 
\vskip.1in
(c) The hypothesis on $s$, $\,s-1/2=k\in\Z^{+}\,$ is not essential, but it does simplify the exposition around the use of the smoothing effect in \cite{KPV1}, which in the inhomogeneous case roughly speaking provides a gain of one derivative, see Theorem \ref{linear1}, estimate \eqref{2.3}.

\vskip.1in

(d) From the assumptions \eqref{newspace} and \eqref{main-2} one has that
$$
\ji x\jd^{m-1}u_0\in  L^2(\R),\;\;\;\;\;\;\;\ji x\jd^{m}u_0\notin  L^2(\R)
$$
and by Sobolev embedding
$$
\ji x\jd^{m}\partial_xu_0,\dots,\ji x\jd^{m}\partial_x^Mu_0\in  L^{\infty}(\R).
$$

(e) We observe that the local smoothing effect \eqref{main-4} obtained in Theorem \ref{main} is slightly weaker than that obtained in \cite{LiPoSa} which was 
$$\,\sup_{x\in\R}\,\int_0^T|\partial_x^{k+1}u(x,t)|^2dt <\infty.
$$
\end{remark}
\vskip.15in

 In \cite{LiPoSa} we proved, roughly, Theorem \ref{main}  with $M=2$ assuming that $\,\nu\leq \epsilon$ for some $\epsilon=\epsilon(\a;\lambda)>0$ small enough.  One of the main tool used in \cite{LiPoSa}  was the homogeneous and inhomogeneous versions of the so called Kato smoothing effect \cite{Ka} for the free Schr\"odinger group $\,\{e^{it\p_x^2}\,:\,t\in\R\}$ (see also \cite{KF}, \cite{CS}, \cite{Ve}, \cite{Sj}) obtained  in \cite{KPV98}. This allows us to apply  the contraction mapping principle  to the integral equation equivalent form for the IVP \eqref{gdnls},
\begin{equation}\label{lower-3}
u(t)=e^{it\p_x^2}u_0+\mu\, \int\limits_0^t e^{i(t-t')\p_x^2} (|u|^{\a}\p_xu)(t')\,dt'.
\end{equation}

These smoothing effects can be expressed as follows : if $\,u=u(x,t)\,$ is the solution of the linear inhomogeneous IVP
\begin{equation}\label{ls}
\begin{cases}
\p_tu=i\,\p_x^2u+f(x,t), \hskip15pt x,t\in\R, \\
u(x,0)=u_0(x),
\end{cases}
\end{equation}
then it satisfies 

\begin{equation}
\begin{aligned}
\label{smooth}
&\sup_{t\in\R}\|u(t)\|_2 \,+\,\sup_{x\in \R}( \int_{-\infty}^{\infty} |D_x^{1/2}u(x,t)|^2dt)^{1/2}\\
\\
&\leq c(\|u_0\|_2\,+\,
\int_{-\infty}^{\infty}\big( \int_{-\infty}^{\infty} |D_x^{-1/2}f(x,t)|^2\,dt\big)^{1/2}dx.
\end{aligned}
\end{equation}

Roughly, since the nonlinearity in the equation in \eqref{gdnls} involves a derivative of order one the inequality \eqref{smooth}
allows to close the estimate avoiding a loss of derivatives. However, in this process a term involving  the 
$\,\|\cdot\|_{L^1(\R:L^{\infty}[0,T])}$-norm appears which can not be made \lq\lq small" by taking $T $ small. So to complete the estimate one needs here  to assume to be working with \lq\lq small" solutions (corresponding to small data).

To overcome this obstruction  we follow the argument given in \cite{KPV98} and rewrite the equation in \eqref{gdnls} as
\begin{equation}\label{eq2}
\p_tu=i\,\p_x^2u +\mu\,|u_0|^{\a}\p_xu +\mu\,(|u|^{\a}-|u_0|^{\a})\p_xu, \hskip10pt x,t\in\R, \hskip5pt 0<\a \le 1,
\end{equation}
and use that under appropriate assumptions the solutions of the linear part of \eqref{eq2} 
\begin{equation}
\label{eq3}
\p_tu=i\,\p_x^2u +\mu\,|u_0|^{\a}\p_xu +f(x,t),
\end{equation}
exhibit smoothing effects similar to those described in \eqref{smooth} (for more detail see Theorem \ref{linear1} in section 2).
These assumptions compromise the regularity and the decay of the lower order coefficient in \eqref{eq3}.
It is here that the number $M$ in the statement of Theorem \ref{main} appears. 

The result in Theorem \ref{main} extends to a related family of equations:

\begin{theorem}\label{thm2}  Under the same hypotheses, the conclusions of Theorem \ref{main} extend to solutions of the IVP associated to the equation
$$
\p_tu=i\p_x^2u + \mu\,\p_x(|u|^{\a}u), \hskip10pt x,t\in\R,  \hskip5pt 0<\a \le 1\;\; {\rm and}\;\; |\mu|=1.
$$
\end{theorem}

This paper is organized as follows, in Section 2, we list some estimates useful in the proof of Theorem \ref{main}. Section 3 contains the proof of our main result Theorem \ref{main}.

Before leaving this section we will introduce the notation to be used in this manuscript.

\subsection{Notation} We denote $\ji x\jd= (1+x^2)^{1/2}$. The Fourier transform of a function $f$, and its inverse Fourier transform are denoted by 
$\hat{f}$ and $\check{f}$ respectively. For $s\in \R$, $ J^s=(1-\p_x^2)^{s/2}$, and $D^s=(-\p_x^2)^{s/2}$ stand for the Riesz and Bessel potentials of
order $-s$, respectively. The functional space $H^s(\R)=(1-\p_x^2)^{-s/2}(L^2(\R))$ denotes the Sobolev spaces of order $s$ endowed with the norm $\|f\|_{s,2}=\|J^s f\|_2$.

For $\,j\in\Z$ we shall use the notation
$$
I_j=[j,j+1],\;\;\;\;\;\;\;\;I_j^T=[j,j+1]\times [0,T],\;\;\;T>0,
$$
and for  two variable functions $f=f(x,t)$ with $(x,t)\in \R\times [0,T]$, 
the norms for $\,1\leq p, q<\infty$ with the usual modification when  $\,p=\infty$ or $\,q=\infty$,
$$
\|f\|_{l^pL^q(I_j^T)}=\Big( \sum_{j\in\Z}\Big( \,\int_{I_j}\,\int_0^T |f(x,t)|^q\,dt dx\,\Big)^{p/q}\Big)^{1/p},
$$
\begin{equation*}
\|f\|_{L^q_TL^p_x}=\Big(\int\limits_0^T\Big(\int\limits_{\R} |f(x,t)|\,dx\Big)^{q/p}\,dt\Big)^{1/q},
\end{equation*}
and
\begin{equation*}
\|f\|_{L^p_xL^q_T}=\Big(\int\limits_{\R}\Big(\int\limits_{0}^T |f(x,t)|\,dt\Big)^{p/q}\,dx\Big)^{1/p}.
\end{equation*}

\section{Linear Estimates}

In this section we shall consider first the IVP for the linear Schr\"odinger equation with lower order variable coefficients. We shall recall some linear estimates obtained in \cite{KPV98} which will be the main tool in the proof of Theorem \ref{main}.To simplify the exposition we shall restrict ourselves to the reduced setting where these estimates are needed, i.e. in one-space dimension ($n=1$) and only one variable coefficient. However, we remark that these assumptions are not essential.
Thus, we consider the IVP
\begin{equation}\label{linearIVP}
\begin{cases}
\p_t u=i\partial_x^2 u+b(x)\, \partial_x u+f_1(x,t)+f_2(x,t),\hskip10pt x\in\R, \;t> 0,\\
u(x,0)=u_0(x),
\end{cases}
\end{equation}
where $\,b:\R\to\C$.
The well-posedness of the IVP \eqref{linearIVP} was studied in several works where necessary and sufficient conditions on the decay and regularity of the coefficient $b(x)$ were deduced. In particular, in the one dimensional case for $\,f_1=f_2\in L^1([0,\infty):L^2(\R))$ with 
$b\in C^1_b(\R)$, i.e. $b\in C^1(\R),\;b, \,b'\in L^{\infty}(\R)$ the condition 
\begin{equation}\label{MC1}
 \underset{l\in \R}{\underset{x\in\R}{\rm Sup}}\;\Big|\int\limits_0^l {\rm Im} \, b(x\pm r\,)\,dr\Big| <\infty
\end{equation}
has been  shown to be  sufficient condition for the $L^2$-well-posedness \cite{Ichinose1}, see also   \cite{Miz}, \cite{Ichinose} and \cite{Takeuchi}.

Our assumptions on the coefficient $\,b(x)\,$ are:
\vskip.1in
(i)  $b\in C^1_b(\R) $ with 
$b:\R \to \C$ satisfying  that
\begin{equation}\label{MC2}
 \underset{l\in \R}{\underset{x\in\R}{\rm Sup}}\;\Big|\int\limits_0^l {\rm Im} \, b(x\pm r\,)\,dr\Big| <\infty.
\end{equation}

(ii) There exists $M\in\Z^{+}$ such that $b\in C^{M}_b(\R)$ with
\begin{equation}
\label{hyp1}
\| b\|_{C^{M}_b} =\sum_{k=0}^M \| b^{(k)}\|_{\infty}=A_1.
\end{equation}
(iii) In addition, 
\begin{equation}
\label{hyp2}
b(x)=\sum_{j\in\Z}\,\alpha_j\,\varphi_j(x),\;\;\text{with}\;\;supp(\varphi_j)\subset [j-1,j+2],\;\;\| \varphi_j\|_{C^{M}}\leq 1,
\end{equation}
with
\begin{equation}
\label{hyp3}
\sum_{j\in\Z}\,|\alpha_j|=A_2.
\end{equation}

The following theorem  is a particular case of a result established  in \cite{KPV98} (Corollary 3.5)
\begin{theorem}\label{linear1}

Under the hypotheses \eqref{MC2}-\eqref{hyp3} for any 
\begin{equation}
\label{2.1}
(u_0,f_1,f_2)\in (H^{1/2}(\R)\times l^1(L^2(I_j^T))\times L^1([0,T]:H^{1/2}(\R)))
\end{equation}
the IVP \eqref{linearIVP} has a unique solution
\begin{equation}
\label{2.2}
u\in C([0,T]:H^{1/2}(\R))
\end{equation}
such that
\begin{equation}
\label{2.3}
\begin{aligned}
&\,\sup_{t\in[0,T]}\|  D_x^{1/2} u(t)\|_2 + \|\,\partial_x u\|_{l^{\infty}L^2(I_j^T)}\\
&\leq c (\| D^{1/2}_xu_0\|_2 +\|f_1\|_{l^{1}L^2(I_j^T)}+\int_0^T\,\|D^{1/2}_xf_2(t)\|_2\,dt),
\end{aligned}
\end{equation}
where $\,c=c(A_1;A_2;T)$.
\end{theorem}

We shall also  work with the following simplified version of Theorem \ref{linear1}:

\begin{corollary}\label{linearA}

Under the hypotheses \eqref{MC2}-\eqref{hyp3} for any 
\begin{equation}
\label{2.4}
(u_0,f_2)\in (L^2(\R)\times L^1([0,T]:L^2(\R)))
\end{equation}
the IVP \eqref{linearIVP} with $\,f_1\equiv 0$ has a unique solution
\begin{equation}
\label{2.5}
u\in C([0,T]:L^2(\R))
\end{equation}
such that
\begin{equation}
\label{2.6}
\sup_{t\in[0,T]}\|  u(t)\|_2 \leq c \,(\|u_0\|_2 +\int_0^T\,\|f_2(t)\|_2\,dt),
\end{equation}
where $\,c=c(A_1;A_2;T)$.
\end{corollary}

Next, we shall introduce some notation. Let $v=v(x,t)$ be the solution of the linear homogeneous IVP 
\begin{equation}\label{linearhomIVP}
\begin{cases}
\p_t v=i\partial_x^2 v+b(x)\, \partial_x v,\hskip10pt x\in\R, \;t\in\R,\\
v(x,0)=v_0(x),
\end{cases}
\end{equation}
with $\,b :\R\to \C$ satisfying the hypotheses \eqref{MC2}-\eqref{hyp3}. 
We define
\begin{equation}
\label{not1}
W_b(t)v_0 = v(\cdot,t),\;\;\;\;t\in\R,
\end{equation}
therefore $\,\{W_b(t)\,:\,t\in\R\}$ is a family of $L^2$-bounded operator (see \eqref{2.6}), with the group property, i.e. 
$$W_b(t) W_b(t')=W_b(t+t'),\hskip5pt t,\,t'\in\R,\hskip10pt \text{and} \hskip10pt W_b(0)=I.$$
Moreover, if $\,v_0\in H^2(\R)$, then
\begin{equation}
\label{der}
\frac{d\:}{dt} W_b(t)v_0= (i\partial_x^2+b(x)\p_x)W_{b}(t) v_0= W_b(t)(i\partial_x^2+b(x)\p_x)v_0.
\end{equation}

In particular, since
\begin{equation}\label{Z1}
\begin{cases}
\p_t (\p_xv)=i\partial_x^2 (\p_xv)+b(x)\, \partial_x (\p_xv) +b'(x)\p_x v,\\
\p_x v(x,0)=\p_xv_0(x),
\end{cases}
\end{equation}
using \eqref{2.6} and for $\,T>0$ it follows that
\begin{equation}
\label{Z21}
\sup_{t\in[0,T]}\|\,\p_xv(t)\|_2\leq c\,( \|\p_x v_0\|_2+\int_0^T \,\|b'\,\p_xv(t)\|_2\,dt).
\end{equation}

Hence, if 
\begin{equation}
\label{Z22}
\, c\,T\,\|b'\|_{\infty}\leq 1/2,
\end{equation}
 one gets that
\begin{equation}
\label{Z3}
\sup_{t\in[0,T]}\|\,\p_xW_b(t)v_0\|_2\leq c \|\p_x v_0\|_2.
\end{equation}
A similar argument shows that if \eqref{Z22} holds, then
\begin{equation}
\label{Z4}
\sup_{t\in[0,T]}\|\,\p_x^2W_b(t)v_0\|_2\leq c\, (\|\p_x^2 v_0\|_2 +T\,\|b''\|_{\infty}\,\|\p_xv_0\|_2) ,
\end{equation}
and
\begin{equation}
\begin{aligned}
\label{Z5}
\sup_{t\in[0,T]}\|\,\p_x^3W_b(t)v_0\|_2&\leq c \,\big(\|\p_x^3 v_0\|_2 +T\,\|b''\|_{\infty}(\|\p_x^2v_0\|_2+T\|b''\|_{\infty}\|\p_xv_0\|_2)\\
&\;\;\;+
T\|b^{(3)}\|_{\infty}\|\p_xv_0\|_2\big).
\end{aligned}
\end{equation}

Formally, the solution of the IVP \eqref{linearIVP} can be written as in the equivalent integral equation form 
\begin{equation}
\label{2.7}
u(t)=W_b(t)u_0 +\int_0^t \,W_b(t-t')(f_1+f_2)(t')\,dt'.
\end{equation}

As it was mentioned in the introduction we shall rewrite the IVP \eqref{gdnls} as
\begin{equation}\label{gdnls2}
\begin{cases}
\p_tu=i\,\p_x^2u +\mu\,|u_0|^{\a}\p_xu +\mu\,(|u|^{\a}-|u_0|^{\a})\p_xu, \;\,  0<\a \le 1,\\
u(x,0)=u_0(x).
\end{cases}
\end{equation}
Therefore, we shall consider the family of operators
$$
\{W_{\mu |u_0|^{\a}}(t)\,:\,t\geq 0\},
$$
for which in order to simplify the notation we shall omit the sub-index $\,\mu\,|u_0|^{\a}$. Thus, we  fix  $b(x)=\mu |u_0|^{\a}$ with $u_0$ as in Theorem \ref{main}.

To obtain weighted estimates for the solution of the IVP 
\begin{equation}\label{abc}
\begin{cases}
\p_tu=i\,\p_x^2u +\mu\,|u_0|^{\a}\p_xu, \:\; 0<\a \le 1,\\
u(x,0)=u_0(x),
\end{cases}
\end{equation}
we multiply the equation in \eqref{abc} by $\,x^n,\;n\in \Z^+$, to get that
\begin{equation}\label{abc1a}
\p_t(x^n u)=i\,\p_x^2(x^n u) +\mu\,|u_0|^{\a}\p_x(x^n u) +h_{0,n}(x,t), 
\end{equation}
with
\begin{equation}
\label{abc2}
h_{0,n}(x,t)= -2 i n x^{n-1}\,\p_xu-in(n-1) x^{n-2} u-\mu\,|u_0|^{\a}\ n x^{n-1} u
\end{equation}

Hence, combining Corollary \ref{linearA} (see \eqref{2.6}) and the hypotheses \eqref{sob}-\eqref{newspace}-\eqref{main-1}  on  $u_0$ in Theorem \ref{main}  it follows that
for any $n\in\Z^+$
\begin{equation}
\label{abc5}
\begin{aligned}
&\sup_{t\in[0,T]}(\| x^{n}u(t)\|_2+\|u(t)\|_2) \\
&\leq c\,\big( \| x^{n}u_0\|_2 +\|u_0\|_2+\int_0^T\| x^{n-1}\p_x u\|_2(t)\,dt\big).
\end{aligned}
\end{equation} 

Next, we see that
\begin{equation}\label{abc6}
\p_t(x^{n-1}\p_x u)=i\,\p_x^2(x^{n-1}\p_x u) +\mu\,|u_0|^{\a}\p_x(x^{n-1}\p_x u) +h_{1,n-1}(x,t), 
\end{equation}
with
\begin{equation}
\label{abc7a}
\begin{aligned}
h_{1,n-1}(x,t)= &-2 i (n-1) x^{n-2}\,\p_x^2u-i(n-1)(n-2) x^{n-3} \p_xu\\
&+\mu\,\p_x(|u_0|^{\a}) x^{n-1} \p_x u-(n-1)x^{n-2}\mu\,|u_0|^{\a} \p_x u.
\end{aligned}
\end{equation}

Thus, by iteration one gets that
\begin{equation}\label{abc8}
\p_t(\p^n_x u)=i\,\p_x^2(\p_x^n u) +\mu\,|u_0|^{\a}\p_x(\p_x^n u) +h_{n,0}(x,t), 
\end{equation}
with
\begin{equation}
\label{abc9}
h_{n,0}(x,t)=\mu\,\sum_{l=1}^nc_r\,\p_x^l(|u_0|^{\alpha})\,\p_x^{n-l+1}u.
\end{equation}

Hence, by combining the interpolation inequalities below (see \eqref{ine1}-\eqref{ine2}) and the hypotheses \eqref{sob}-\eqref{main-2} we deduce that there exists $\,T=T(\lambda, \nu,\alpha)>0$ such that for any $\,j\in\Z^+,\;j=1,\dots,m$
\begin{equation}
\label{abc10}
\sup_{t\in[0,T]}(\| \,x^j\,u(t)\|_2+\|u(t)\|_{j,2})\leq c\,(\| x^j u_0\|_2+\| u_0\|_{j,2}).
\end{equation}

\vskip.1in
As a consequence of Corollary \ref{linearA} and the interpolation estimates in \eqref{ine1}-\eqref{ine2} 
using the hypotheses of Theorem \ref{main} one obtains that the solution of the IVP
\begin{equation}
\label{A1}
\begin{cases}
\p_t v=i\partial_x^2 v+\mu |u_0|^{\alpha} \partial_x v+f(x,t),\hskip10pt x\in\R, \;t\in[0,T],\;\;T>0,\\
v(x,0)=v_0(x),
\end{cases}
\end{equation}
satisfies:

\begin{lemma}\label{linear3}
Let $m, k, M$ be as in Theorem \ref{main}. Then the following estimates hold for the solution of the IVP \eqref{A1} 
\begin{equation}
\label{a24}
\sup_{t\in[0,T]} \tres v(t)\tres \leq c(\tres v_0\tres +\int_0^T\,\tres f(t)\tres\,dt),
\end{equation}
where
\begin{equation}
\label{a3}
\tres v(t)\tres \equiv\, \sum_{j=1}^{M+1}\| |x|^m\,\partial_x^j v(t)\|_2\,+\,\sum_{j=0}^{s-1/2}\| \partial_x^j v(t)\|_2.
\end{equation}
\end{lemma}

The proof of Lemma \ref{linear3} follows an argument similar to that delineated in \eqref{abc}-\eqref{abc10}.

\vskip.1in
To end this section,  we state the following useful interpolation results. 
\begin{lemma}\label{interp} For any $a, b>0$ and $\gamma\in (0,1)$ 
\begin{equation}
\begin{aligned}
\label{ine1}
&\|J^{\gamma a}(\ji x\jd^{(1-\gamma)b}f)\|_2 \le  c\,\|\ji x\jd^{b} f\|_2^{1-\gamma}\|J^a f\|^{\gamma}_2, \hskip10pt \gamma\in(0,1),\\
\text{and\hskip15pt}&\\
&\|\ji x\jd^{\gamma a}(J^{(1-\gamma)b}f)\|_2 \le  c\,\|J^{b} f\|_2^{1-\gamma}\|\ji x\jd^a f\|^{\gamma}_2, \hskip10pt \gamma\in(0,1).
\end{aligned}
\end{equation} 
\end{lemma}

For the proof which is based on the Three Lines Theorem we refer to \cite{NP}.

By integration by parts, one also has :
\begin{lemma}\label{interp2} For any $j, \,k\in\Z^{+},\,j,\,k\geq 1$ there exists $c=c(k;j)>0$ such that
\begin{equation}
\begin{aligned}
\label{ine2} 
& \| \ji x\jd^k\partial_x^jf\|_2^2\leq c \| \ji x\jd^k\partial_x^{j+1}f\|_2 \| \ji x\jd^k\partial_x^{j-1}f\|_2 +c \| \ji x\jd^{k-1}\partial_x^{j-1}f\|_2^2,\\
& \| \ji x\jd^k\partial_x^jf\|_2^2\leq c \| \ji x\jd^{k-1}\partial_x^{j+1}f\|_2 \| \ji x\jd^{k+1}\partial_x^{j-1}f\|_2 +c \| \ji x\jd^{k-1}\partial_x^{j-1}f\|_2^2,\\
& \| \ji x\jd^k\partial_x^jf\|_2^2\leq c \| \ji x\jd^{k+1}\partial_x^{j+1}f\|_2 \| \ji x\jd^{k-1}\partial_x^{j-1}f\|_2 +c \| \ji x\jd^{k-1}\partial_x^{j-1}f\|_2^2.
\end{aligned}
\end{equation} 
\end{lemma}
\vskip.15in

\section{Proof of Theorem \ref{main}}

To solve the IVP \eqref{gdnls} we consider its  integral version form
\begin{equation}\label{a003}
u(t)= W(t) u_0+\int_0^t W(t-t') (|u|^{\alpha}-|u_0|^{\alpha})\p_x u (t')\,dt
\end{equation}
where $W(t)$ is the solution of the linear problem
\begin{equation}\label{a004}
\begin{cases}
\p_t w= i\p_x^2w + \mu |u_0|^{\alpha}\p_x w\\
w(x,0)=w_0(x).
\end{cases}
\end{equation}

Thus, we shall prove that the operator
\begin{equation}
\label{w1}
\Phi (v(t)) = W(t) u_0+\int_0^t W(t-t') (|v|^{\alpha}-|u_0|^{\alpha})\p_xv (t')\,dt,
\end{equation}
defines a contraction in the space
\begin{equation}
\label{espace}
\begin{split}
X_T&=\Big\{ v\in C([0,T]: H^s(\R)) : \\
&\underset{[0,T]}{\sup}\;\big(\|v(t)\|_{s,2}+\| \ji x\jd^m v(t)\|_{\infty}+\sum_{j=0}^M\| \ji x\jd^m \partial_x^{j+1}v(t)\|_2\big) \\
\\
&\;\;\;+ \|\p_x^{k+1}v\|_{\ell^{\infty}L^2(I_j^T)}\le 2c\nu,\\
\\
&\underset{[0,T]}{\sup}  \,\sum_{j=0}^{1}\| \ji x\jd^m \partial_t\p_x^{j}v(t)\|_2\leq 2c\nu(1 +(2c\nu)^{2}),\\
\\
\;\;\;\;\;&\text{and}\\
&\;\;\sup_{[0,T]}\|  \ji x\jd^m (v(\cdot,t)-u_0)\|_{\infty}\leq  \frac{\lambda}{4}\;\;\Big\}
\end{split}
\end{equation}
for $\,T>0\,$ sufficiently small. 

We observe that if
\begin{equation}
\label{qq1}
\sup_{[0,T]}\|  \ji x\jd^m (v(\cdot,t)-u_0)\|_{\infty}\leq  \frac{\lambda}{4},
\end{equation}
then from the hypothesis \eqref{main-2} one has that
\begin{equation}
\label{qq2}
\inf_{(x,t)\in\R\times[0,T]}\;\ji x\jd^m |v(x,t)|\geq \frac{\lambda}{2}.
\end{equation}
Also since $\Phi(v(t))$ denotes the solution of the linear IVP
\begin{equation}
\label{71}
\begin{cases}
\begin{aligned}
&\p_t w=i\p_x^2 w +\mu|u_0|^{\alpha} \p_xw+\mu(|v|^{\alpha}-|u_0|^{\alpha})\p_xv,\\
&w(x,0)=u_0(x),
\end{aligned}
\end{cases}
\end{equation}
the estimates for $\,\underset{[0,T]}{\sup}\| \ji x\jd^m \partial_t\p_x^{j}v(t)\|_2,\;j=0,1\,$ will follow directly from those not involving derivatives in the $t$-variable. 

Hence, we need to estimate $\Phi$ in the following norms in \eqref{espace}
\begin{equation}\label{a005}
\tres v \tres + \sup_{[0,T]}\|\p_x^kD_x^{1/2}v\|_2+ \|\p_x^{k+1}v\|_{\ell^{\infty}L^2(I_j^T)}+\sup_{[0,T]} \|\ji x\jd^m v\|_{\infty},
\end{equation}
with $ \,\tres \cdot \tres$  defined in \eqref{a3}. 
The inequality \eqref{a24} provides the estimate
\begin{equation}\label{a005b}
\begin{split}
\tres \Phi(v)\tres  &\lesssim \tres u_0\tres +\int_0^T\tres (|v|^{\alpha}-|u_0|^{\alpha})\p_x v\tres\,dt\\
&  \lesssim \tres u_0\tres + T\sup_{[0,T]} \tres |v|^{\alpha}\p_xv\tres+ T\sup_{[0,T]} \tres |u_0|^{\alpha}\p_xv\tres.
\end{split}
\end{equation} 
The last two terms can be handled using the argument presented in detail in \cite{LiPoSa}.
Thus it  remains to estimate $\Phi(v) $ in the last three norms in \eqref{a005}.

From \eqref{2.3}  one has that
\begin{equation}\label{a007}
\begin{split}
\underset{[0,T]}{\sup} \|\p_x^{k+1/2}\Phi(v)\|_2&+ \|\p_x^{k+1}\Phi(v)\|_{\ell^{\infty}L^2(I_j^T)} \\
&\lesssim \|u_0\|_{s,2}
+\overset{k-1}{\underset{l=0}{\sum}} \|\p_x^l(|v|^{\alpha}-|u_0|^{\alpha})\p_x^{k-l}v\|_{\ell^1L^2(I_j^T)}\\
&\hskip10pt+\|\p_x^k(|v|^{\alpha}-|u_0|^{\alpha})\p_x v\|_{\ell^1L^2(I_j^T)}\\
&\hskip10pt+ \|(|v|^{\alpha}-|u_0|^{\alpha})\p_x^{k+1}v\|_{\ell^1L^2(I_j^T)}.
\end{split}
\end{equation}

Next we will restrict to analyze the last two terms since the other ones are easier to deal with (or follows similar arguments).

An application of Holder's inequality yields
\begin{equation}\label{a008}
 \begin{split}
& \|(|v|^{\alpha}-|u_0|^{\alpha})\p_x^{k+1}v\|_{\ell^1L^2(I_j^T)}\\
&\lesssim \Big\| \int_0^t \frac{d}{d\theta}|v|^{\alpha}\,d\theta\Big\|_{\ell^1 L^{\infty}(I_j^T)}
 \|\p_x^{k+1} v\|_{\ell^{\infty}L^2(I_j^T)}\\
 &\lesssim \Big\|\int_0^T \big  |v|^{\alpha-1}\p_t v\, d\theta\Big\|_{\ell^1 L^{\infty}(I_j^T)} \|\p_x^{k+1} v\|_{\ell^{\infty}L^2(I_j^T)}\\
 &\lesssim  T \sup_{[0,T]} \| |v|^{\alpha-1}\p_t v\,\ji x\jd \|_{1,2} \|\p_x^{k+1} u\|_{\ell^{\infty}L^2(I_j^T)}\\
 &\lesssim T \big(\| \ji x\jd^m\p_t v\|_{L^{\infty}_TL^2} +\| \ji x\jd^m\p_t\p_x v\|_{L^{\infty}_TL^2} \\
 &\hskip15pt+\|\ji x\jd^m\p_x^2v\|_{L^{\infty}_TL^2_x}\big)\|\p_x^{k+1} v\|_{\ell^{\infty}L^2(I_j^T)},
 \end{split}
 \end{equation}
 

 where we have used that
 \begin{equation}\label{a009}
 \sum_{j} a_j \le c \Big(\sum_{j} \ji j\jd^2 a_j^2\Big)^{1/2}
 \end{equation}
 and Sobolev embedding to obtain
 \begin{equation}\label{a010}
 \sum_{j} \sup_{I_j}|f|\lesssim \sum_{j} \|f\|_{H^1(I_j)}\lesssim \Big(\sum_{j}\ji j\jd ^2\|f\|_{H^1(I_j)}^2\Big)^{1/2}\simeq \|\ji x\jd f\|_{H^1(\R)}.
 \end{equation}
 Also above we have utilized  that
 \begin{equation}
 \label{abc1}
 \inf_{\R\times[0,T]}\ji x\jd^m |v(x,t)|\geq \frac{\lambda}{2}
 \end{equation}
 to estimate the negative powers of $\,|v(x,t)|$.
 
 To handle $\|\p_x^k(|v|^{\alpha}-|u_0|^{\alpha})\p_x v\|_{\ell^1L^2(I_j^T)}$ we restrict ourselves to the terms that might present more difficulties
 \begin{equation}\label{a011}
\|\big(|v|^{\alpha-1}\p_x^kv-|u_0|^{\alpha-1}\p_x^ku_0)\, \p_x v\|_{\ell^1L^2(I_j^T)}
\end{equation}
 and
\begin{equation}\label{a012}
\|\big(|v|^{\alpha-k}(\p_x v)^k-|u_0|^{\alpha-k}(\p_xu_0)^k)\, \p_x v\|_{\ell^1L^2(I_j^T)}.
 \end{equation}
 
Employing \eqref{a009} it follows that
\begin{equation}\label{a013}
\begin{split}
&\|\big(|v|^{\alpha-1}\p_x^kv-|u_0|^{\alpha-1}\p_x^ku_0)\, \p_x v\|_{\ell^1L^2(I_j^T)}\\
&\le \||v|^{\alpha-1}\p_x^kv\, \p_x v\|_{\ell^1L^2(I_j^T)}+\||u_0|^{\alpha-1}\p_x^ku_0\, \p_x v\|_{\ell^1L^2(I_j^T)}\\
&\lesssim  \|\ji x\jd |v|^{\alpha-1}\p_x^kv\, \p_x v\|_{L^2_TL^2_x}+\|\ji x\jd |u_0|^{\alpha-1}\p_x^ku_0\, \p_x v\|_{L^2_TL^2_x}\\
&\lesssim \|\ji x\jd^{(1-\alpha)m +1}\p_x^kv\, \p_x v\|_{L^2_TL^2_x}+\|\ji x\jd^{(1-\alpha)m +1}\p_x^ku_0\, \p_x v\|_{L^2_TL^2_x}\\
&\lesssim T^{1/2} \big(\|\ji x\jd^m \p_x v\|_{L^{\infty}_TL^{\infty}_x}\|\p_x^k v\|_{L^{\infty}_TL^2_x}
+\|\ji x\jd^m \p_x v\|_{L^{\infty}_TL^{\infty}_x}\|\p_x^k u_0\|_{L^{\infty}_TL^2_x}\big).
\end{split}
\end{equation}

\medskip

As in the previous estimate we use \eqref{a009} to lead to
\begin{equation}\label{a014}
\begin{split}
&\|\big(|v|^{\alpha-k}(\p_x v)^k-|u_0|^{\alpha-k}(\p_xu_0)^k)\, \p_x v\|_{\ell^1L^2(I_j^T)}\\
&\lesssim  \||v|^{\alpha-k}(\p_x v)^k\p_x v\|_{\ell^1L^2(I_j^T)} +\||u_0|^{\alpha-k}(\p_xu_0)^k\, \p_x v\|_{\ell^1L^2(I_j^T)}\\
&\lesssim  \|\ji x\jd^{(k-\alpha)m+1}(\p_x v)^k\p_x v\|_{L^2_TL^2_x} +\|\ji x\jd^{(k-\alpha)m+1}(\p_xu_0)^k\, \p_x v\|_{L^2_TL^2_x}\\
&\lesssim  T^{1/2} \big(\|\ji x\jd^m \p_xv\|_{L^{\infty}_TL^{\infty}_x}^k
+\|\ji x\jd^m \p_x u_0\|_{L^{\infty}_TL^{\infty}_x}^k\big)\|\p_x v\|_{L^{\infty}_TL^2_x}.
\end{split}
\end{equation}

Combining \eqref{a007}, \eqref{a008}, \eqref{a013}, \eqref{a014} and interpolation we lead to
\begin{equation}\label{a014a}
\begin{split}
\underset{[0,T]}{\sup} \|\p_x^{k+1/2}&\Phi(v)\|_2+ \|\p_x^{k+1}\Phi(v)\|_{\ell^{\infty}L^2(I_j^T)} \\
& \lesssim \|u_0\|_{s,2}+ c T\nu^2((1+(2c\nu)^2)+c\,T^{1/2}\nu^2+ cT^{1/2}\nu^{k+1}.
\end{split}
\end{equation}

\vskip.1in

Next we need to take care of the norm $\| \ji x\jd^m\Phi(v)\|_{\infty}$ and to show that  
\begin{equation}\label{main-2ab}
\sup_{[0,T]} \|\ji x\jd^m(\Phi(v)-u_0)\|_{\infty}\le \frac{\lambda}{4}.
\end{equation}

We recall the hypothesis \eqref{main-2}
\begin{equation}\label{main-2a}
\underset{x}{\rm Inf}\; \ji x\jd^m\,|u_0(x)|\ge \lambda>0.
\end{equation}
Hence, formally one has that for any $\,t\in [0,T]$ with $\,T$ such that 
\begin{equation}
\label{Q1}
T\,\| \p_x|u_0|^{\alpha}\|_{\infty}\leq c \,T \,\| \ji x\jd^{(1-\alpha)m}\,\p_xu_0\|_{\infty}\leq 1/2,
\end{equation}
(see \eqref{Z22}) it follows that
\begin{equation}\label{AA1}
\begin{aligned}
\|W(t)u_0-u_0\|_{\infty}&=\|\,\int_0^t \frac{d}{d\tau} \,W(\tau)u_0\,d\tau\|_{\infty}\\
&=\|\,\int_0^t  \,W(\tau)(i\p_x^2+\mu|u_0|^{\alpha} \p_x)u_0\,d\tau\|_{\infty}\\
&\leq \|\,\int_0^t  \,W(\tau)(i\p_x^2+\mu|u_0|^{\alpha} \p_x)u_0\,d\tau\|_{1,2}\\
&\leq  c\,T\,\|u_0\|_{3,2},
\end{aligned}
\end{equation}
after combining Sobolev embedding and the estimates \eqref{2.6} and \eqref{Z22}

Thus, assuming \eqref{Q1} for any $t\in[0,T]$ with $T$ sufficiently small and following the argument leading to
\eqref{abc10} we get
\begin{equation}\label{AA2}
\begin{aligned}
& \|\ji x\jd^m \,(W(t)u_0-u_0)\|_{\infty} = \|\ji x\jd^m \int_0^t  \,\frac{d}{d\tau}W(\tau)u_0\,d\tau\|_{\infty}\\
&\leq  \|\ji x\jd^m \int_0^t  \,W(\tau)(i\p_x^2+\mu|u_0|^{\alpha} \p_x)u_0\,d\tau\|_{\infty}\\
&\leq c T \sup_{[0,T]} \|\ji x\jd^m  \,W(t)(i\p_x^2+\mu|u_0|^{\alpha} \p_x)u_0\|_{\infty}\\
&\leq cT\,(\, \sum_{j=1}^3\| \ji x\jd^m \p_x^j u_0\|_2+\|u_0\|_{m+3,2})\leq cT \,\nu\leq \frac{\lambda}{4}.
\end{aligned}
\end{equation}

\medskip

Therefore, for $\,T\,$ sufficiently small

\begin{equation}
\label{pqr1}
\begin{aligned}
\sup_{[0,T]}\|\ji x\jd^m \,(\Phi(v(t))&-u_0)\|_{\infty} \leq \sup_{[0,T]}\|\ji x\jd^m \,(W(t)u_0-u_0)\|_{\infty}\\
&+ \sup_{[0,T]}\|\ji x\jd^m \int_0^tW(t-t')(|v|^{\alpha}-|u_0|^{\alpha})\p_x v(t')dt'\|_{\infty}\\
&\leq \;\;\;\frac{\lambda}{4} + cT (c\nu)(1+(c\,\nu))\leq \frac{\lambda}{2}.
\end{aligned}
\end{equation}

Hence, 
\begin{equation}
\label{abc7}
\inf_{\R\times[0,T]}\ji x\jd^m|\Phi(v)(x,t)|\geq \frac{\lambda}{2},
\end{equation}
see \eqref{qq1}-\eqref{qq2}.

\medskip

From the argument above we also have that
\begin{equation}\label{abc7b}
\|\ji x\jd^m\Phi(v)\|_{\infty}\le  c\|\ji x \jd^m u_0\|_{\infty}+ cT(c\nu) (1+(c\nu)).
\end{equation}

Therefore,  combining \eqref{abc7b}, \eqref{a005b}, \eqref{a014a} and the remark in \eqref{71} we deduce
\begin{equation}\label{abc7c}
\begin{split}
\tres \Phi(v)\tres &+ \|\ji x\jd^m\Phi(v)\|_{L^{\infty}_TL^{\infty}_x}\\
&\hskip10pt +\underset{[0,T]}{\sup}  \,\sum_{j=0}^{1}\| \ji x\jd^m \partial_t\p_x^{j}\Phi(v)\|_2
+ \|\p_x^{k+1}\Phi(v)\|_{\ell^{\infty}L^2(I_j^T)} \\
&\le c\nu + c T^{1/2}( \nu^2 +\nu^{k+1})+ cT\nu(1+(2c\nu)^2)+cT(c\nu) (1+(c\nu))\\
&\le 2c\nu
\end{split}
\end{equation}
whenever
\begin{equation}\label{Texis}
c T^{1/2}( \nu +\nu^{k})+ cT(1+(2c\nu)^2)+cT (1+(c\nu)) \ll 1.
\end{equation}

This shows that the map $\Phi$ is well defined.
The argument employed above also allows us to show that $\Phi(v)$ is a contraction. Hence we will omit the details.

\medskip

We notice that the continuous dependence of the initial data  does not follows from the contraction mapping principle argument we just
have employed above.

Next we shall show the continuous dependence of the initial data. To do so, we consider $u, v$ solutions of the IVP with initial data  $u_0$ and $v_0$  respectively, i.e.
\begin{equation}\label{a1}
\p_t u = i\, \p_x^2u +\mu\, |u_0|^{\alpha}\,\p_x u +\mu\, (|u|^{\alpha}-|u_0|^{\alpha})\, \p_x u.
\end{equation}
\begin{equation}\label{a21}
\p_t v = i\, \p_x^2v +\mu\, |v_0|^{\alpha}\,\p_x v +\mu\, (|v|^{\alpha}-|v_0|^{\alpha})\, \p_x v.
\end{equation}

Consider now $w=u-v$, then $w$ satisfies the equation
\begin{equation}\label{a33}
\begin{split}
\p_t w=& i\, \p_x^2w +\mu\, |u_0|^{\alpha}\,\p_x w +\mu\, (|u_0|^{\alpha}-|v_0|^{\alpha})\, \p_x v \\
&+ \mu(|u|^{\alpha}-|u_0|^{\alpha}- (|v|^{\alpha}-|v_0|^{\alpha}))\, \p_x u+ \mu  (|v|^{\alpha}-|v_0|^{\alpha})\, \p_x w.
\end{split}
\end{equation}

 After differentiate $k$ times with respect to $x$,  the more delicate terms are
 \begin{equation}\label{a4}
 \|\big(|u|^{\alpha}-|u_0|^{\alpha}- (|v|^{\alpha}-|v_0|^{\alpha}))\, \p_x^{k+1} u\|_{\ell^1L^2(I_j^T)}\equiv A_1
 \end{equation}
 and
  \begin{equation}\label{a5}
 \|\p^k_x\big(|u|^{\alpha}-|u_0|^{\alpha}- (|v|^{\alpha}-|v_0|^{\alpha}))\, \p_x u\|_{\ell^1L^2(I_j^T)}\equiv A_2
 \end{equation}
 
 We start by estimating $A_1$. Using Holder's inequality
 \begin{equation}\label{a6}
 \begin{split}
 A_1&\lesssim \Big\| \int_0^t \frac{d}{d\theta}|u|^{\alpha}\,d\theta- \int_0^t \frac{d}{d\theta}|v|^{\alpha}\,d\theta\Big\|_{\ell^1 L^{\infty}(I_j^T)}
 \|\p_x^{k+1} u\|_{\ell^{\infty}L^2(I_j^T)}\\
 &\lesssim \Big\|\int_0^T \big | |u|^{\alpha-1}\p_t u-|v|^{\alpha-1}\p_tv\big|\,d\theta\Big\|_{\ell^1 L^{\infty}(I_j^T)} \|\p_x^{k+1} u\|_{\ell^{\infty}L^2(I_j^T)}\\
 &\lesssim  T \sup_{[0,T]} \| (|u|^{\alpha-1}\p_t u-|v|^{\alpha-1}\p_tv)\ji x\jd \|_{1,2} \|\p_x^{k+1} u\|_{\ell^{\infty}L^2(I_j^T)}
 \end{split}
 \end{equation}
 where we have used the inequalities \eqref{a009} and \eqref{a010}.

 Next we estimate the first norm on the right hand side of the last inequality in \eqref{a6}.
 
 \begin{equation}\label{a6b}
 \begin{split}
\underset{[0,T]}{\sup} \| & \ji x\jd (|u|^{\alpha-1}\p_t u-|v|^{\alpha-1}\p_tv))\|_{1,2}\\
&\lesssim  \| \ji x\jd (|u|^{\alpha-1}\p_tu-|v|^{\alpha-1})\p_tv \|_{L^{\infty}_T L^2_x}\\
&\hskip10pt+ \| \ji x\jd \p_x(|u|^{\alpha-1}\p_tu-|v|^{\alpha-1}\p_tv) \|_{L^{\infty}_T L^2_x}\\
&= A_{1,1}+A_{1,2}.
\end{split}
\end{equation}

We start by estimating $A_{1,1}$.

Recalling that
 \begin{equation}\label{a14}
 |u|^{\alpha-1}-|v|^{\alpha-1}= c\, \big( \theta |u|+(1-\theta)|v|\big)^{\alpha-2} \big| |u|-|v|\big|
 \end{equation}
 and
 \begin{equation}\label{a15}
 \ji x\jd^m\,|u| \ge \lambda,\;  \ji x\jd^m\,|v| \ge \lambda,
 \end{equation}
we have 
\begin{equation}\label{a15b}
\begin{split}
A_{1,1} &\lesssim \| \ji x\jd (|u|^{\alpha-1}-|v|^{\alpha-1}\p_tu\|_{L^{\infty}_TL^2_x} +\|\ji x\jd |u|^{\alpha-1}\p_t(u-v)\|_{L^{\infty}_TL^2_x}\\
&\lesssim \| \ji x\jd^{m(2-\alpha)+1} |u-v| \p_tu\|_{L^{\infty}_TL^2_x} +\|\ji x\jd^{m(1-\alpha)+1}\p_t w\|_{L^{\infty}_TL^2_x} \\
&\lesssim \| \ji x\jd^m w\|_{L^{\infty}_TL^{\infty}_x} \| \ji x\jd^m \p_tu\|_{L^{\infty}_TL^2_x} +\|\ji x\jd^m \p_t w\|_{L^{\infty}_TL^2_x}.
\end{split}
\end{equation}

On the other hand
\begin{equation}\label{a15c}
\begin{split}
A_{1,2} &\le \| \ji x\jd (|u|^{\alpha-2}\p_xu\p_tu - |v|^{\alpha-2}\p_xv\p_t v\|_{L^{\infty}_TL^2_x}\\
&\hskip10pt + \|\ji x\jd (|u|^{\alpha-1}\p_t\p_xu-|v|^{\alpha-2}\p_t\p_x v\|_{L^{\infty}_TL^2_x}\\
&= A_{1,2,1} + A_{1,2,2}.
\end{split}
\end{equation}
 
 A  similar argument as the one used in \eqref{a15b} yields
 \begin{equation}\label{a15d}
 \begin{split}
 &A_{1,2,2} \le \|\ji x\jd (|u|^{\alpha-2}-|v|_{\alpha-2})\p_xu\p_t u\|_{L^{\infty}_TL^2_x}\\
 &\hskip10pt+\|\ji x\jd |v|^{\alpha-2}\p_xu\p_t(u-v)\|_{L^{\infty}_TL^2_x}+\|\ji x\jd \p_x(u-v)\p_tv \|_{L^{\infty}_TL^2_x}\\
 &\lesssim \|\ji x\jd^{m(3-\alpha)+1}|u-v| \p_xu\p_tv\|_{L^{\infty}_TL^2_x}\\
 &\hskip10pt \| \ji x\jd ^{m(2-\alpha)+1}\p_xu\p_tw\|_{L^{\infty}_TL^2_x}+\|\ji x\jd^{m(2-\alpha)+1}\p_xw\p_tv\|_{L^{\infty}_TL^2_x}\\
 &\lesssim  \|\ji x\jd^m w\|_{L^{\infty}_TL^{\infty}_x} \|\ji x\jd^m \p_xu\|_{L^{\infty}_TL^{\infty}_x} \|\ji x\jd^m\p_tv\|_{L^{\infty}_TL^2_x}\\
 &\hskip10pt + \|\ji x\jd^m \p_x u\|_{L^{\infty}_TL^{\infty}_x}  \|\ji x\jd^m\p_t w\|_{L^{\infty}_TL^2_x}
 + \|\ji x\jd^m \p_x w\|_{L^{\infty}_TL^{\infty}_x} \|\ji x\jd^m\p_tv\|_{L^{\infty}_TL^2_x}.
 \end{split}
 \end{equation}
 
 We employ a similar argument to get
 \begin{equation}\label{a15e}
 A_{1,2,2} \lesssim  \|\ji x\jd^m w\|_{L^{\infty}_TL^{\infty}_x} \|\ji x\jd^m \p_x\p_tu\|_{L^{\infty}_TL^2_x} + \| \ji x\jd^m\p_t\p_x w\|_{L^{\infty}_TL^2_x}.
 \end{equation}

 \smallskip
 
 Now we proceed to estimate $A_2$, we first see that
 \begin{equation}\label{a9}
 \begin{split}
 A_2 &\le \|\p^k_x\big(|u|^{\alpha}-|v|^{\alpha})\, \p_x u\|_{\ell^1L^2(I_j^T)}+ \|\p^k_x\big(|u_0|^{\alpha}-|v_0|^{\alpha})\, \p_x u\|_{\ell^1L^2(I_j^T)}.\\
 \end{split}
 \end{equation}
 We observe that the last term contains the expression $\p^k_x\big(|u_0|^{\alpha}-|v_0|^{\alpha})$ which does not depend on time and so the analysis will
 not be so difficult. Thus we will be focus on the first term on the right hand side.
 
 As before we estimate the two more difficult terms in $\|\p^k_x\big(|u|^{\alpha}-|v|^{\alpha})\, \p_x u\|_{\ell^1L^2(I_j^T)} $, that is,
\begin{equation}\label{a10}
\|\big(|u|^{\alpha-1}\p_x^ku-|v|^{\alpha-1}\p_x^kv)\, \p_x u\|_{\ell^1L^2(I_j^T)}\equiv A_{2,1}
\end{equation}
 and
\begin{equation}\label{a11}
\|\big(|u|^{\alpha-k}(\p_x u)^k-|v|^{\alpha-k}(\p_xv)^k)\, \p_x u\|_{\ell^1L^2(I_j^T)}\equiv A_{2,2}.
 \end{equation}
 
 We have that
 \begin{equation}\label{a12}
 A_{2,1}\lesssim \||u|^{\alpha-1}\p_x^k(u-v)\p_xu\|_{\ell^1L^2(I_j^T)}+\|(|u|^{\alpha-1}-|v|^{\alpha-1})\p_x^kv\, \p_xu\|_{\ell^1L^2(I_j^T)}.
 \end{equation}
 
 Using \eqref{a008}  and \eqref{a15} we deduce that
 \begin{equation}\label{a13}
 \begin{split}
\||u|^{\alpha-1}\p_x^k(u-v)\p_xu\|_{\ell^1L^2(I_j^T)}&\lesssim \|\ji x\jd |u|^{\alpha-1}\p_x^k w \,\p_xu\|_{L^2_TL^2_x}\\
&\lesssim \|\ji x\jd \ji x\jd^{(1-\alpha)m}\p_x u \,\p_x^k w\|_{L^2_TL^2_x}\\
&\lesssim T^{1/2} \|\ji x\jd^m \p_x u\|_{L^{\infty}_TL^{\infty}_x}\| \p_x^k w\|_{L^{\infty}_TL^2_x}.
\end{split}
\end{equation}

Using \eqref{a14} and \eqref{a15}, and applying \eqref{a008} we obtain 
 \begin{equation}\label{a16}
 \begin{split}
\|(&|u|^{\alpha-1}-|v|^{\alpha-1})\p_x^kv\, \p_xu\|_{\ell^1L^2(I_j^T)}\\
&\lesssim \| \ji x\jd^{m(2-\alpha)+1} |u-v| \p_x^k v \,\p_xu\|_{L^2_TL^2_x}\\
&\lesssim  T^{1/2} \|\ji x\jd^m \p_x u\|_{L^{\infty}_TL^{\infty}_x}\|\p_x^kv\|_{L^{\infty}_TL^2_x}\|\ji x\jd^m w\|_{L^{\infty}_TL^{\infty}_x}.
\end{split}
\end{equation} 
 
 Finally,  $A_{2,2}$ can be estimated as follows: we first use \eqref{a008} to obtain
 \begin{equation}\label{a17}
 \begin{split}
 A_{22} &\lesssim \|\ji x\jd \Big[ (|u|^{\alpha-k}-|v|^{\alpha-k})(\p_xu)^k+|v|^{\alpha-k}\big( (\p_xu)^k-(\p_xv)^k\big)\Big]\p_xu\|_{L^2_TL^2_x}\\
 &\lesssim \|\ji x\jd (|u|^{\alpha-k}-|v|^{\alpha-k})(\p_xu)^k\p_xu\|_{L^2_TL^2_x}\\
 &\hskip10pt + \| \ji x\jd|v|^{\alpha-k}\big( (\p_xu)^k-(\p_xv)^k\big)\p_xu\|_{L^2_TL^2_x}.
 \end{split}
 \end{equation}
 Then applying \eqref{a14} and \eqref{a15} once more we obtain
 \begin{equation}\label{a18}
 \begin{split}
 \|\ji x\jd (|u|^{\alpha-k}-&|v|^{\alpha-k})(\p_xu)^k\p_xu\|_{L^2_TL^2_x} \\
 &\lesssim \|\ji x\jd^{m(k+1-\alpha)+1}|w|(\p_xu)^{k+1}\|_{L^2_TL^2_x}\\
 &\lesssim T^{1/2} \| \ji x\jd^m \p_x u\|_{L^{\infty}_TL^{\infty}_x}^{k+1}\|w\|_{L^{\infty}_TL^2_x},
 \end{split}
 \end{equation}
 and
 \begin{equation}\label{a19}
 \begin{split}
 &\| \ji x\jd|v|^{\alpha-k}\big( (\p_xu)^k-(\p_xv)^k\big)\p_xu\|_{L^2_TL^2_x}\\
 &= c \| \ji x\jd|v|^{\alpha-k}\, \p_x w \big[ (\p_xu)^{k-1}+\cdots+(\p_xv)^{k-1}\big]\p_xu\|_{L^2_TL^2_x} \\
 &\lesssim T^{1/2} \|\p_x w\|_{L^{\infty}_TL^2_x} \|\ji x\jd^m\p_xu\|_{L^{\infty}_TL^{\infty}_x}
  \sum_{l=0}^{k-1} \|\ji x\jd^m\p_xu\|_{L^{\infty}_TL^{\infty}_x}^{k-1-l}\|\ji x\jd^m\p_xv\|_{L^{\infty}_TL^{\infty}_x}^l.
 \end{split}
 \end{equation}

 In summary, combining  the estimates  \eqref{a6} to \eqref{a19} we have basically proved
 \begin{equation}\label{19b}
 \begin{split}
 &\tres w\tres +\underset{[0, {T'}]}{\sup} \|\p_x^{k+1/2}w\|_2+ \|\p_x^{k+1}w\|_{\ell^{\infty}L^2(I_j^{T'})}\\
 &\hskip10pt + \| \ji x \jd^m w\|_{L^{\infty}_{T'}L^{\infty}_x}+\underset{[0,T']}{\sup}  \,\sum_{j=0}^{1}\| \ji x\jd^m \partial_t\p_x^{j}w\|_2\\
 & \le c \tres u_0-v_0\tres+\|u_0-v_0\|_{s,2}
 \end{split}
 \end{equation}
 for $T'\in (0,T)$  with $T$ given in \eqref{Texis} which implies our claim.



\vspace{.5cm} 
\begin{thebibliography}{99}



\bibitem{Linares} H. Biagioni and F. Linares, {\it Ill-posedness for the derivative Schr\"odinger and generalized Benjamin-Ono equations},
Trans. Amer. Math. Soc. {\bf 353} (2001), 3649--3659.


\bibitem{Cazenave} T. Cazenave and I. Naumkin, {\it Local existence, global existence, and scattering for the nonlinear Schr\"odinger equation},
 Comm. Contemp. Math. {\bf 19} (2017),  1650038, 20 pp.


\bibitem{Colin} M. Colin and  M. Ohta, {\it Stability of solitary waves for derivative nonlinear Schr\"odinger equation},  Ann. Inst. H. Poincar\'e Anal. Non Lin\'eaire {\bf 23} (2006), 753--764.  

\bibitem{Iteam} J. Colliander, M. Keel, G. Staffilani, H. Takaoka and T. Tao, {\it Global well-posedness for the Schr\"odinger equations with derivative},
SIAM J. Math. Anal. {\bf 33} (2001), 649--669. 

\bibitem{refinement} J. Colliander, M. Keel, G. Staffilani, H. Takaoka and T. Tao, {\it A refined global well-posedness for the Schr\"odinger equations with derivative}, SIAM J. Math. Anal. {\bf 34} (2002), 64--86.

\bibitem{CS}  P. Constantin and J.-C. Saut,  {\it Local smoothing properties of dispersive equations}, J. Amer. Math. Soc. {\bf 1} (1988), 413--446.


\bibitem{FHI} N. Fukaya, M. Hayashi, and T. Inui, {\it A sufficient condition for global existence of solutions to a generalized derivative nonlinear Schr\"odinger equation}, Anal. PDE {\bf 10} (2017), 1149--1167.


\bibitem{Grunrock} A. Gr\"unrock and S. Herr, {\it Low regularity local well-posedness of the derivative nonlinear Schr\"odinger equation with periodic initial data},  SIAM J. Math. Anal.  {\bf 39} (2008), 1890--1920.  



\bibitem{Guo} B. Guo and  Y. Wu, {\it Orbital stability of solitary waves for the nonlinear derivative Schr\"odinger equation},  J. Differential Equations  {\bf 123} (1995), 35--55.


\bibitem{GW} Z. Guo and Y. Wu, {\it Global well-posedness for the derivative nonlinear Schr\"odinger equation in $H^{1/2}(\R)$}, Discrete Cont. Dyn. Syst.
{\bf 37} (2017),  257--264.

\bibitem{Hao} C. Hao, {\it Well-posedness for one-dimensional derivative nonlinear Schr\"odinger equations}, Comm. Pure Appl. Anal. {\bf 6} (2007), 997--1021.





 
 

\bibitem{HayashiII} N. Hayashi, {\it The initial value problem for the derivative nonlinear  Schr\"odinger equation in the energy space}, Nonlinear Anal.  TMA  {\bf 20} (1993), 823--833.

 \bibitem{MHO} M. Hayashi and T. Ozawa, {\it Well-posedness for a generalized derivative nonlinear Schr\"odinger equation},  J. Differential Equations
  {\bf 261} (2016), 5424--5445.

\bibitem{HayashiOzawa} N. Hayashi and T. Ozawa, {\it On the derivative nonlinear Schr\"odinger equation}, Phys. D {\bf 55} (1992), 14--36.

\bibitem{HayashiOzawa2} N. Hayashi and T. Ozawa, {\it Finite energy solution of nonlinear Schr\"odinger  equations of derivative type}, 
SIAM J. Math. Anal. {\bf 25} (1994), 1488--1503.


\bibitem{Herr} S. Herr, {\it On the Cauchy problem for the derivative nonlinear Schr\"odinger equation with periodic boundary condition}, Int. Math. Res. Not. IMRN (2006), Art. ID 96763, 33 pp.

\bibitem{Ichinose} W. Ichinose, {\it On the Cauchy problem for Schr\"odinger type equations and Fourier integral operators}, J. Math. Kyoto Univ. {\bf 33} (1993), , 583--620

\bibitem{Ichinose1} W. Ichinose, {\it On a necessary condition for $L^2$-well-posedness of the Cauchy problem for some Schr\"odinger type equations with a potential term}, J. Math. Kyoto Univ. {\bf 33} (1993),  647--663


\bibitem{JLPS} R. Jenkins, J. Liu, P. Perry, and C. Sulem, {\it Global well-posedness and soliton resolution for the derivative nonlinear Schr\"odinger equation},  arXiv: 1706.06252v1.


\bibitem{Ka} T. Kato, {\it On the Cauchy problem for the (generalized) Korteweg-de Vries equation}, 
  Advances in Mathematics Supplementary Studies, Studies in Applied Math. {\bf 8} (1983), 93--128


 \bibitem{KN} D. J. Kaup and A. C. Newell, {\it An exact solution  for a derivative nonlinear Schr\"odinger equation}. J. Math. Phys.  {\bf 19} (1978), 789--801.

 \bibitem{KPV1} C. E. Kenig, G. Ponce, L. Vega, {\it Small solutions to nonlinear Schr\"odinger equations}, Ann. Inst. H. Poincar\'e Anal. 
 Non Lin\'eaire {\bf 10} (1993), 255--288.

\bibitem{KPV98} C. E. Kenig, G. Ponce, L. Vega, {\it Smoothing effects and local existence theory for the generalized nonlinear Schrödinger equations}, Invent. Math. 134 (1998), no. 3, 489?545.

 \bibitem{KF} S. N. Kruzhkov and A. V. Faminskii, {\it Generalized solutions of the Cauchy problem for the Korteweg-de Vries
    equation}, Math. U.S.S.R. Sbornik  {\bf 48} (1984), 93--138.

\bibitem{KW} S. Kwon and Y. Wu, {\it Orbital stability of solitary waves for derivative nonlinear Schr\"odinger equation}, preprint, arXiv:1603.03745.

\bibitem{LW} S. Le Coz and  Y. Wu, {\it Stability of multi-solitons for the derivative nonlinear Schr\"odinger equation}, preprint, arXiv:1609.04589.

\bibitem{LP} F. Linares and G. Ponce, {\it Introduction to Nonlinear Dispersive Equations} (Second Edition). Universitext. Springer, New
York, 2015.

\bibitem{LiPoSa} F. Linares, G. Ponce, and G. N. Santos, {\it On a class of solutions to the generalized derivatives Schr\"odinger equations}, to appear in Acta Math. Sinica English series.

\bibitem{LSS} X. Liu, G. Simpson and C. Sulem, {\it Stability of solitary waves for a generalized derivative nonlinear Schr\"odinger equation}, 
 J. Nonlinear Sci. {\bf 23} (2013), no. 4, 557--583.
 
 \bibitem{MTX} C. Miao, X. Tang, and  G.  Xu, {\it Stability of the traveling waves for the derivative Schr\"odinger equation in the energy space}, preprint
  arXiv:1702.07856.
  
  
 \bibitem{Mio} K. Mio, T. Ogino, K. Minami and S. Takeda, {\it Modified nonlinear Schr\"odinger equation for Alfv\'en Waves propagating along magnetic field in cold plasma}, J. Phys. Soc. {\bf 41} (1976), 265--271.

\bibitem{Miz} S. Mizohata, {\em  On the Cauchy problem}, Notes and Reports in Mathematics in Science and Engineering {\bf 3}, Academic Press
Inc., Orlando, FL, 1985.


\bibitem{Mijolhus} E. Mjolhus,  {\it On the modulational instability of hydromagnetic waves parallel to the magnetic field}, J. Plasma Phys. {\bf 16} (1976), 321--334.

\bibitem{MMW} J. Moses, B. A. Malomed, F. W. Wise, {\it Self-steepening of ultrashort optical pulses without self-phase-modulation}, Phys. Rev. A {\bf 76} (2007), 1--4.


\bibitem{NP} J. Nahas and G. Ponce, {\it On the Persistent Properties of Solutions to Semi-Linear Schr\"odinger Equation}, Comm. Partial Differential Equations {\bf 34}  (2009), 1208--1227.

\bibitem{Ozawa} T. Ozawa, {\it On the nonlinear Schr\"odinger equations of derivative type}, Indiana Univ. Math. J. {\bf 45} (1996), 137--163.


\bibitem{Gleison}  G. N.  Santos, {\it Existence and uniqueness of solutions for a generalized Nonlinear Derivative Schr\"odinger equation},  J. Differential Equations  {\bf 259} (2015), 2030--2060. 


\bibitem{Sj}  P. Sj\"olin, {\it Regularity of solutions to the Schr\"odinger equations}, Duke Math. J.  {\bf 55} (1987), 699--715

\bibitem{Takaoka} H. Takaoka, {\it Well-posedness for the one-dimensional nonlinear Schr\"odinger equation with the derivative nonlinearity}, Adv. Differential Equations  {\bf 4} (1999), 561--580.

\bibitem{Takaoka2} H. Takaoka, {\it Global well-posedness for the Schr\"odinger equations with derivative in a nonlinear term and data in low-order Sobolev spaces},   Electron. J. Differential Equations {\bf 43} (2001), 1--23.


\bibitem{Takeuchi} J. Takeuchi, {\it  A necessary condition for the well-posedness of the Cauchy problem for a certain class of evolution equations},   Proc. Japan Acad.  {\bf 50} (1974), 133--137.

\bibitem{TX} X. Tang and G.  Xu, {\it Stability of the sum of two solitary waves for (gDNLS) in the energy space}, preprint  arXiv:1702.07858. 

\bibitem{Tsutsumi} M. Tsutsumi and I. Fukuda, {\it On solutions of the derivative nonlinear Schr\"odinger equation. Existence and uniqueness theorem}, Funkcialaj Ekvacioj {\bf 23} (1980), 259--277.


\bibitem{Ve} L. Vega, {\it The Schr\"odinger equation: pointwise convergence to the initial data}, 
 Proc. Amer. Math. Soc. {\bf 102} (1988), 874--878.

\bibitem{W} Y. Wu, {\it Global well-posedness of the derivative nonlinear Schr\"odinger equation}, Analysis \& PDE {\bf 8} (2015), 1101--1112. 

\end{thebibliography}
\end{document}